\documentclass[11pt]{amsart}
\usepackage{fullpage,array,bbm, graphicx, amssymb, bbm, rotating,epic,multirow,pst-node, mathrsfs, stmaryrd}
\usepackage{amsmath}
\usepackage{amsthm}

 \font\tenmsb=msbm10 at 12pt \font\sevenmsb=msbm7 at 8pt \font\fivemsb=msbm5 at
6pt
\newfam\msbfam
\textfont\msbfam=\tenmsb \scriptfont\msbfam=\sevenmsb \scriptscriptfont\msbfam=\fivemsb

\def\R{{\mathbb R}}

\def \F{{\mathscr R}}
\def \T{{\mathscr T}}

\def\N{{\mathbb N}}

\def\Z{{\mathbb Z}}

\begin{document}
\newcommand{\reset}{\setcounter{equation}{0}}

\newcommand{\beq}{\begin{equation}}
\newcommand{\noi}{\noindent}
\newcommand{\eeq}{\end{equation}}
\newcommand{\dis}{\displaystyle}
\newcommand{\mint}{-\!\!\!\!\!\!\int}

\def \theequation{\arabic{section}.\arabic{equation}}

\newtheorem{thm}{Theorem}[section]
\newtheorem{lem}[thm]{Lemma}
\newtheorem{cor}[thm]{Corollary}
\newtheorem{prop}[thm]{Proposition}
\theoremstyle{definition}
\newtheorem{defn}{Definition}
\newtheorem{rem}[thm]{Remark}

\def \bx{\hspace{2.5mm}\rule{2.5mm}{2.5mm}} \def \vs{\vspace*{0.2cm}} \def
\hs{\hspace*{0.6cm}}
\def \ds{\displaystyle}
\def \p{\partial}
\def \O{\Omega}
\def \b{\beta}
\def \m{\mu}
\def \ou{{\overline u}}
\def \ov{{\overline v}}
\def \D{\Delta}
\def \d{\delta}
\def \s{\sigma}
\def \e{\varepsilon}
\def \a{\alpha}
\def \o{\omega}
\def \g{\gamma}
\def \wv{{\widetilde v}}
\def \hw{{\widehat w}}
\def\cqfd{%
\mbox{ }%
\nolinebreak%
\hfill%
\rule{2mm} {2mm}%
\medbreak%
\par%
}
\def \pr {\noindent {\it Proof:} }
\def \rmk {\noindent {\it Remark} }
\def \esp {\hspace{4mm}}
\def \dsp {\hspace{2mm}}
\def \ssp {\hspace{1mm}}

\title{One dimensional sharp discrete Hardy-Rellich inequalities}
\author{Xia Huang}
\address{School of Mathematical Sciences,  Key Laboratory of MEA (Ministry of Education) \& Shanghai Key Laboratory of PMMP,  East China Normal University, Shanghai 200241, China}
\email{xhuang@cpde.ecnu.edu.cn}
\author{Dong Ye}
\address{School of Mathematical Sciences,  Key Laboratory of MEA (Ministry of Education) \& Shanghai Key Laboratory of PMMP,  East China Normal University, Shanghai 200241, China}
\address{IECL, UMR 7502, University of Lorraine, 57050 Metz, France}
\email{dye@math.ecnu.edu.cn}

\date{}
\begin{abstract}
In this paper, we establish discrete Hardy-Rellich inequalities on $\mathbb{N}$ with $\Delta^\frac{\ell}{2}$ and optimal constants, for any $\ell \geq 1$. As far as we are aware, these sharp inequalities are new for $\ell \geq 3$. Our approach is to use weighted equalities to get some sharp Hardy inequalities using shifting weights, then to settle the higher order cases by iteration. We provide also a new Hardy-Leray type inequality on $\N$ with the same constant as the continuous setting. Furthermore, the main ideas work also for general graphs or the $\ell^p$ setting.
\end{abstract}

\maketitle
\vskip 0.8cm

\section{Introduction}

\reset

In 1921, in a letter to Hardy, Landau \cite{Lau} gave a proof of the following inequality with optimal constants: For any $p > 1$ and any non-negative real sequence ${\{a_n}\}_{n \geq 1}$, we have
\begin{equation}\label{H1}
 \left(\frac{p-1}{p}\right)^p \sum_{n=1}^\infty\left(\frac{a_1+a_2+...+a_n}{n}\right)^p \leq \sum_{n=1}^\infty a_n^p.
\end{equation}

\medskip
Let $C_c(\N)$ be the space of finitely supported real sequences indexed by the set of nonnegative integers $\N$, or equivalently, the space of finitely supported real functions on $\N$. With abuse of notations, we will write $u_n = u(n) = (u)_n$ for any sequence $u$ over $\N$. We denote also the discrete gradient over $\N$ by $\nabla u_n := (\nabla u)_n  = u_n - u_{n-1}$ for $n \geq 1$, then \eqref{H1} is equivalent to
\begin{equation}\label{p}
\sum_{n=1}^\infty |\nabla u_n|^p \geq \left(\frac{p-1}{p}\right)^p \sum_{n=1}^\infty \frac{|u_n|^p}{n^p}, \quad \forall\; u\in C_c(\N),\; u_0=0.
\end{equation}
The most well-known case with $p=2$, was established by Hardy \cite{H} when he proposed a simple proof for a Hilbert's inequality.
Since one century, Hardy type inequalities have been enriched extensively and broadly, they play an important role in many branches of analysis and geometry. A huge literature exists, especially for the continuous setting, we refer to the classical or more recent books \cite{BEL, GM, M, OK, RS} for interested readers.

\medskip
Although the history of Hardy inequalities found its origin somehow in the discrete setting, the discrete Hardy-Rellich inequalities are much less understood comparing to the continuous situation.
In 2018, Keller-Pinchover-Pogorzelski \cite{KPP0, KPP} proved the following improved discrete Hardy inequality of  first order:
\begin{equation}\label{kpp1}
\sum_{n=1}^\infty |\nabla u_n|^2 \geq  \sum_{n=1}^\infty w_n u_n^2, \quad \forall\; u\in C_c(\N),\; u_0=0
\end{equation}
with $w_n = \frac{\Delta \sqrt{n}}{\sqrt{n}}$, where $\Delta v_n= 2v_n - v_{n+1} - v_{n-1}$ denotes the discrete Laplacian over $\N$. More explicitly, for $n \geq 1$,
$$w_n = \sum_{k=1}^\infty \binom{4k}{2k}\frac{1}{2^{4k-1}(4k-1)n^{2k}} = \frac{1}{4n^2} + \frac{5}{64n^4}+ \frac{21}{512n^6}+\ldots.$$

Inspired by \cite{DFP}, they showed that the weights $(w_n)$ are critical in many aspects, in particular the inequality \eqref{kpp1} cannot hold with any sequence $\gneq w_n$. See also \cite{FKP} for similar improvement of \eqref{p} with general $p > 1$.

\medskip
 A natural question is to handle the higher order cases, i.e. the Hardy-Rellich inequalities on $\N$ with $\Delta^\frac{\ell}{2}$, $\ell \geq 2$. For odd integers $\ell \geq 1$, by $\Delta^\frac{\ell}{2}$ we mean $\nabla \Delta^\frac{\ell-1}{2}$.

\medskip
Recently, Gupta established in \cite{Gup1} the following Hardy-Rellich inequalities of general order over $C_c(\N)$: Let $\ell \geq 2$,
\begin{align}
\label{Gup-l}
\sum_{2n \geq \ell} |\Delta^\frac{\ell}{2}u_n|^2  \geq 2^{\ell -3}(\ell - 1)! \sum_{n= \ell}^\infty \frac{u_n^2}{n^{2\ell}}, \quad \forall\; u\in C_c(\N), \; u_k = 0, \; 0 \leq k \leq \ell-1.
\end{align}
In \cite[Theorem 2.5]{Gup1}, the LHS of \eqref{Gup-l} was written as the sum on $n \geq 1$, we can see that under the assumptions $u_k = 0$ for $k \leq \ell-1$, $\Delta^\frac{\ell}{2}u_n = 0$ if $2n < \ell$.

\medskip
However, the question of finding the best constants are widely open for $\ell \geq 2$. Very recently, the second order case was solved by Gerhat-Krej\v{c}i\v{r}\'\i k-\v{S}tampach \cite{GKS}, they proved that for $\ell = 2$, the constant in \eqref{Gup-l} can be raised to $\frac{9}{16}$, the sharp one expected for $\Delta$. More precisely, there holds
\begin{align}\label{gks}
\sum_{n\geq 1} |\Delta u_n|^2 \geq \sum_{n\geq 2} \widetilde w_n u^2_n, \quad \forall\; u\in C_c(\N),\; u_0=u_1 =0,
\end{align}
with
$$\widetilde w_n = 6\sum_{k=1}^\infty \frac{4^k - 1}{4^{2k}}\frac{(4k)!}{(2k)!(2k+2)!}\frac{1}{n^{2k+2}} = \frac{9}{16n^4} + \frac{210}{256n^6} + \ldots$$
Their main idea is to use a factorization technique, but this approach seems to be difficultly adaptable for higher order situations.

\medskip
Our main concern here is to provide the following Hardy-Rellich inequalities with optimal constants for any order $\ell \geq 2$.
\begin{thm}
\label{thm1.1}
Let $\ell \geq 2$ be integer. Given any $u\in C_c(\N)$ with $u_k=0$ for $0 \leq k \leq \ell-1$, there holds
\begin{align}
\label{rellich-l}
\sum_{2n \geq \ell } |\Delta^\frac{\ell}{2}u_n|^2  \geq \Big[\frac{(2\ell)!}{4^\ell\ell!}\Big]^2\sum_{n= \ell}^\infty \frac{u_n^2}{n^{2\ell}}.
\end{align}
\end{thm}

Notice that the constant in \eqref{rellich-l} coincides the sharp constant for Hardy-Rellich inequality in the continuum (see for instance \cite{OW}):
\begin{align}
\label{rellich-cl}
\int_0^\infty |\varphi^{(\ell)}(x)|^2 dx \geq \Big[\frac{(2\ell)!}{4^\ell\ell!}\Big]^2\int_0^\infty \frac{\varphi^2(x)}{x^{2\ell}} dx, \quad \forall\; \varphi \in C_c^\ell(0, \infty).
\end{align}
By establishing the link between discrete and continuous Hardy-Rellich inequalities, we will show the optimality of the constant in \eqref{rellich-l}, see Proposition \ref{sharp} below. We notice also that the assumptions $u_k=0$, $0 \leq k \leq \ell-1$ are crucially necessary to expect any Hardy-Rellich inequalities like \eqref{rellich-l}, see Remark \ref{zeros} below, this answers some question in \cite{Gup1}. 

\medskip
Furthermore, our approach can enable us a sequence of weights which are larger than that in \eqref{rellich-l}, as in Keller-Pinchover-Pogorzelski or Gerhat-Krej\v{c}i\v{r}\'\i k-\v{S}tampach's work with $\ell = 1$ and $2$ respectively. In particular, for $\ell=2$, for any $u\in C_c(\N)$ with $u_0=u_1=0$, we obtain that
\begin{align*}
\sum_{n \geq 1} |\Delta u_n|^2  \geq \sum_{n\geq 2} \left(\frac{9}{16n^4} + \frac{15}{16 n^5}+ \frac{213}{128 n^6}+ \cdot\cdot\cdot\right) u^2_n.
\end{align*}
This improves asymptotically the weights given in \cite{GKS}, see section 7.2.

\medskip
We will make use of the following elementary observation: As $\nabla u_n = u_n - u_{n-1}$, the corresponding divergence operator in $C_c(\N)$ is just ${\rm div}(u)_n = u_{n+1} - u_n$, in other words
$${\rm div} = \tau_1\circ \nabla = \nabla\circ \tau_1,$$ where $\tau_k$ states for the shift operator $(\tau_kg)_n= g_{n+k}$ for $k \in \Z$. Readily $\Delta = - {\rm div}\circ \nabla$, so the Laplacian over $\N$ is a combination of shift and gradient operators. Remark also that the three operators $\Delta$, $\nabla$ and ${\rm div}$ over $\N$ commute with each other. Therefore, to get higher order Rellich inequalities, we need just to well understand the first order situation with suitable weights. Our main ideas are getting Hardy inequalities via equalities, combined with iteration and shift.

\medskip
In fact, we establish the following first order Hardy inequalities with shifting weights.
\begin{thm}
\label{thm1.2}
Let $\alpha \in (-\infty, 0)$. Then
\begin{align}
\label{a<0bis}
\sum_{n=1}^{\infty} n^\alpha |\nabla v_n|^2 \geq \frac{(\alpha-1)^2}{4} \sum_{n=1}^{\infty} (n+1)^{\alpha-2}v_n^2, \quad \forall\; v\in C_c(\N), \; v_0=0,
\end{align}
or equivalently
\begin{align}
\label{a<0}
\sum_{n=2}^{\infty} (n-1)^\alpha |\nabla u_n|^2 \geq \frac{(\alpha-1)^2}{4} \sum_{n=2}^{\infty} n^{\alpha-2}u_n^2, \quad \forall\; u\in C_c(\N), \; u_0=u_1=0.
\end{align}
\end{thm}

\medskip
Our departure point is the following general equality, which can be called {\sl weighted ground state transform}. Let $f_n>0$ for $n\geq 1$ and $f_0=0$, then for any weights $V = (V_n)$ and any $u\in C_c(\N)$ with $u_0=0$, there holds
\begin{align}
\label{equa1}
\sum_{n=1}^\infty V_n|\nabla u_n|^2 + \sum_{n=1}^\infty \frac{{\rm div}(V\nabla f)_n}{f_n} u^2_n =\sum_{n=2}^\infty V_n\left|\sqrt{\frac{f_{n-1}}{f_n}}u_n - \sqrt{\frac{f_n}{f_{n-1}}}u_{n-1}\right|^2.
\end{align}
Indeed, let $n\geq 2$, direct calculation yields
\begin{align*}
\left|\sqrt{\frac{f_{n-1}}{f_n}} u_n -\sqrt{\frac{f_n}{f_{n-1}}} u_{n-1}\right|^2
& = (u_n - u_{n-1})^2 -  \frac{ f_n-f_{n-1}}{f_n} u_n^2 +  \frac{f_n-f_{n-1}}{f_{n-1}}u_{n-1}^2\\
& = |\nabla u_n|^2 - \frac{\nabla f_n}{f_n} u_n^2 + \frac{\nabla f_n}{f_{n-1}} u_{n-1}^2.
\end{align*}
So we get
\begin{align*}
\sum_{n=2}^{\infty}V_n \left|\sqrt{\frac{f_{n-1}}{f_n}} u_n -\sqrt{\frac{f_n}{f_{n-1}}} u_{n-1}\right|^2
& = \sum_{n=2}^{\infty}V_n |\nabla u_n|^2 - \sum_{n=2}^{\infty}V_n \frac{\nabla f_n}{f_n} u_n^2 + \sum_{n=1}^{\infty}V_{n+1}\frac{\nabla f_{n+1}}{f_n} u_n^2\\
& = \sum_{n=2}^{\infty} V_n |\nabla u_n|^2 + \sum_{n=2}^{\infty}\frac{{\rm div}(V\nabla f)_n}{f_n} u_n^2 + V_2\frac {\nabla f_2}{f_1} u_1^2.
\end{align*}
As $u_0=f_0=0$, then $\nabla u_1= u_1$ and ${\rm div}(V_1 \nabla f_1)= V_2 \nabla f_2 - V_1 \nabla f_1,$ this means that \eqref{equa1} holds.

\medskip
The special case of \eqref{equa1} with $V_n \equiv 1$ can be found in \cite{KLS, KS}. With $V_n\equiv 1$ and $f_n =\sqrt n$ for $n\in \N$, as the RHS of \eqref{equa1} is nonnegative, we get readily the Hardy inequality
\begin{equation*}
\sum_{n=1}^\infty |\nabla u_n|^2 \geq  \sum_{n=1}^\infty \frac{\Delta \sqrt{n}}{\sqrt{n}} u_n^2, \quad \forall\; u\in C_c(\N), u_0=0.
\end{equation*}
which is just the optimal estimate \eqref{kpp1} found by Keller-Pinchover-Pogorzelski \cite{KPP}.

\medskip
As interesting application of \eqref{equa1}, we obtain also the following sharp Hardy inequalities.
\begin{thm}
\label{prop1.2'}
Let $\alpha \geq 0$. Then
\begin{align}
\label{H}
\sum_{n=1}^{\infty} n^\alpha |\nabla u_n|^2 \geq \frac{(\alpha-1)^2}{4} \sum_{n=1}^{\infty}  n^{\alpha-2} u_n^2, \quad \forall\; u\in C_c(\N), \; u_0=0.
\end{align}
\end{thm}

To study the Hardy-Rellich inequalities over $\N$, Gupta considered weighted Hardy inequalities by diverse methods. Firstly, based on the ground state transform, he showed in \cite{Gup} that \eqref{H} is valid for $\alpha\in[0,1)\cup[5,\infty)$, and pointed out that this approach fails to work if $\alpha< 0$ or $\alpha\in (1,4)$, see \cite[Remark 1.1]{Gup}; then he proved \eqref{H} for $\alpha \in 2\N$ by Fourier analysis in \cite{Gup1}. The estimate \eqref{H} with $\alpha\in[0,1)$ follows also from \cite[Theorem 1.1]{Cop} with $\kappa=2$ and $c=2-\alpha$. Furthermore, when $\alpha\in[\frac{1}{3}, 1)$, Gupta gave an improved version by adding infinitely many positive lower order terms in RHS.

\medskip
Here we present a unified approach to handle \eqref{a<0} and \eqref{H} for all $\alpha \in \R$. It is worthy to mention that we do not need \eqref{H} to consider higher order Rellich inequalities \eqref{rellich-l}.

\medskip
Clearly, the estimate \eqref{H} is trivial with $\alpha=1$. In that case, we show a new Leray-type inequality over $\N$ with an explicit constant. We can find similar result in \cite[Lemma 4.5]{KL} with unknown constant.
\begin{thm}
\label{leray}
For any $u\in C_c(\N)$ with $u_0=u_1=0$, there holds
\begin{align}
\label{leray1}
\sum_{n=2}^\infty n |\nabla u_n|^2 \geq \frac{1}{4}\sum_{n=2}^\infty \frac{|u_n|^2}{n(\ln n)^2}.
\end{align}
\end{thm}

In \eqref{leray1}, we find again the sharp constant in the continuum, see \cite[Theorem 1]{ww}.
\begin{align*}
\int_1^\infty x|\varphi'(x)|^2 dx \geq \frac{1}{4}\int_1^\infty \frac{\varphi^2(x)}{x(\ln x)^2} dx, \quad \forall\; \varphi \in C_c^1(1, \infty).
\end{align*}

\medskip
Indeed, our approach works for locally finite graphs. In section 2, we briefly describe some notions and provide some general higher order identities on locally finite graphs, see Theorems \ref{thm2.1}-\ref{thm2.2}. The key inequalities \eqref{a<0} and Theorem \ref{thm1.1} will be established in section 3 and 4 respectively. The proof of \eqref{H} is given in section 5. Finally, we use similar ideas to consider the lattices $\mathbb{Z}^d$ ($d \geq 2$) and the $\ell^p$ setting ($p > 1$).

\section{Weighted Hardy-Rellich equalities on graphs}
\reset
Here we recall some standard notions for locally finite graphs $X$, i.e. each vertex $x\in X$ has finite neighbors. A graph over a countable discrete set $X$ is a symmetric function $b: X \times X\to \mathbb{R}_+$ which vanishes at the diagonal.

\medskip
We refer to the elements of $X$ as vertices. We say that $x,y \in X$ are adjacent or neighbors if $b(x,y)>0$, in which case we write $x\sim y$. Assume that $X$ is equipped with the discrete topology and $b$ is a graph over $X$. Moreover, $b$ is said connected, if for every $x$ and $y$ in $X$ there are $x_0,...,x_m$ in $X$ such that $x_0=x$, $x_m=y$ and $x_i \sim x_{i+1}$ for $i=0,...,m-1$. Given a graph $b$ over $X$, the associated Laplacian $\Delta$ is defined as
$$
\Delta u (x):=\sum_{y\sim x} b(x,y)\big[u(x)-u(y)\big].
$$

We will use the formal inner product (when it makes sense) for $f, g \in C(X)$.
$$
\langle f, g\rangle:= \sum_{x\in X}f(x)g(x).
$$
Let $\nabla f(x,y):=f(x)-f(y)$ be the discrete gradient defined on $X \times X$, and define
$$
\langle \nabla f, \nabla g\rangle:=\sum_{x\in X} \langle\nabla f,\nabla g \rangle_x \quad \mbox{with }\; \langle\nabla f,\nabla g \rangle_x :=\frac{1}{2}\sum_{y\sim x}b(x,y)\big[f(x)-f(y)\big] \big[g(x)-g(y)\big].
$$
Denote also
$$
|\nabla f|^2(x) := \langle\nabla f,\nabla f \rangle_x=\frac{1}{2}\sum_{y\sim x} b(x,y)\big[f(x)-f(y)\big]^2.
$$
By direct calculation, we get some elementary equalities. Let $f \in C(X)$, $g \in C_c(X)$, we have
\begin{align*}
\langle \Delta f, g\rangle& = \sum_{x,y\in X} b(x,y)\big[f(x)-f(y)\big]g(x)\\
& = \sum_{x,y\in X} b(x,y)\big[f(x)-f(y)\big] \big[g(x)-g(y)\big] +\sum_{x,y\in X} b(x,y)\big[f(x)-f(y)\big]g(y)\\
& = 2\langle \nabla f, \nabla g\rangle - \langle \Delta f, g\rangle.
\end{align*}
Hence the following Green's formula and Leibniz rule hold true: For any $f \in C(X)$, $g \in C_c(X)$,
\begin{align}
\label{Green}
\langle\Delta f,g\rangle=\langle f, \Delta g\rangle= \langle \nabla f, \nabla g\rangle,
\end{align}
\begin{align}
\label{Leib}
\Delta(fg)(x)= \big[f\Delta g + g\Delta f\big](x) -2\langle\nabla f,\nabla g \rangle_x.
\end{align}
Consequently, for any $h \in C(X)$, $u \in C_c(X)$, we have
\begin{align}
\label{equa23}
\langle \Delta u, u\rangle=\sum_{x\in X}|\nabla u|^2(x) \quad \mbox{and}\quad \langle \Delta h, u^2\rangle=\langle h, \Delta(u^2)\rangle=2\langle h, u\Delta u- |\nabla u|^2\rangle.
\end{align}

Let $f \in C(X)$ be positive, $V \in C(X)$ and $u \in C_c(X)$, we define
\begin{align*}
\F(V, f, u) := \frac{1}{2}\sum_{x, y\in X} b(x, y)V(x)\left(\sqrt{\frac{f(x)}{f(y)}}u(y)-\sqrt{\frac{f(y)}{f(x)}}u(x)\right)^2.
\end{align*}
Based on above notations, there holds
\begin{align*}
\F(V, f, u)=&\frac{1}{2}\sum_{x\in X}V(x)\left[2|\nabla u|^2(x)+\sum_{y\sim x}b(x,y)\frac{f(x)-f(y)}{f(y)}u^2(y)-\frac{\Delta f}{f}(x)u^2(x)\right]\\
=&\sum_{x\in X}V(x)|\nabla u|^2(x) + \sum_{y\sim x}\frac{\langle\nabla f,\nabla V \rangle_y}{f(y)} u^2(y)-\frac{1}{2}\sum_{y\sim x}V(y)\frac{\Delta f}{f}(y)u^2(y)\\&-\frac{1}{2}\sum_{x\in X}V(x)\frac{\Delta f}{f}(x)u^2(x)\\
 = &\sum_{x\in X}\left[V(x)|\nabla u|^2(x) +\frac{\langle\nabla f,\nabla V \rangle_x}{f(x)} u^2(x) - \frac{V\Delta f}{f}(x)u^2(x)\right]\\
 = &\sum_{x\in X} V(x)|\nabla u|^2(x) - \sum_{x\in X} \frac{\T(V, f)}{f}(x) u^2(x)
\end{align*}
where
$$\T(V, f)(x): = V(x)\Delta f(x) - \langle\nabla f,\nabla V \rangle_x, \quad \forall\; x \in X.$$ Hence we obtain a first order Hardy type equality
\begin{align}\label{FH}
\sum_{x\in X} V(x)|\nabla u|^2(x) - \sum_{x\in X} \frac{\T(V, f)}{f}(x) u^2(x) = \F(V, f, u).
\end{align}

Let us define a divergence operator for $F\in C(X\times X)$ as follows.
\begin{align}\label{div}
{\rm div}F(x):= \frac{1}{2}\sum_{y\sim x} b(x,y)\big[F(y,x)-F(x,y)\big].
\end{align}
For $F,~G \in C(X\times X)$, we denote the formal inner product (when it makes sense) as follows.
$$\langle F, G\rangle := \frac{1}{2}\sum_{x, y \in X} b(x,y)F(x,y)G(x, y).$$
By the symmetry of $b$, we check readily that
\begin{align*}
-\langle {\rm div}F, h\rangle =\langle F, \nabla h\rangle, \quad \forall\; F \in C(X\times X), \;h \in C_c(X).
\end{align*}
An interesting fact is
\begin{lem}
\label{div-T}
For any $V, f \in C(X)$ and any $x \in X$, $\T(V, f)(x) = -{\rm div}(V\nabla f)(x)$. Here $$(V\nabla f)(x,y):=V(x)\nabla f(x,y).$$
\end{lem}

\medskip
\noindent
{\bf Proof}. Let $x \in X$, there holds
\begin{align*}
V(x)\Delta f(x)- \langle \nabla V, \nabla f\rangle_x & = V(x)\sum_{y\sim x} b(x,y)\nabla f(x,y) - \frac{1}{2} \sum_{y\sim x} b(x,y) \nabla V(x,y) \nabla f(x,y)\\
& = \frac{1}{2} \sum_{y\sim x} b(x,y)[V(x)+V(y)] \nabla f(x,y)\\
& = \frac{1}{2} \sum_{y\sim x} b(x,y)\big[V(x)\nabla f(x,y) - V(y)\nabla f(y,x)\big]\\
& = -{\rm div}(V\nabla f)(x),
\end{align*}
where the last equality comes from the definition \eqref{div}. \qed

\begin{rem}
As far as we are aware, we did not find a reference of divergence for general function $F \in C(X\times X)$. Clearly, when $F$ is anti-symmetric as the gradient, it coincides with the definition of divergence for discrete differential forms over graphs, see for example \cite{HJ}.
\end{rem}

If $V \geq 0$, we get immediately from \eqref{FH} a weighted Hardy inequality since $\F(V, f, u) \geq 0$.
\begin{align}
\label{FHdiv}
\sum_{x\in X} V(x)|\nabla u|^2(x) \geq \sum_{x\in X} \frac{-{\rm div}(V\nabla f)}{f}(x) u^2(x), \quad \forall\; u \in C_c(X).
\end{align}

\medskip
Let $h\in C(X)$. Applying \eqref{Green}, \eqref{equa23},  \eqref{FH} and Lemma \ref{div-T} with $Vh$, we get for any $u \in C_c(X)$,
\begin{align*}
\sum_{x\in X}V|\Delta u - h u|^2(x) & = \sum_{x\in X} V|\Delta u|^2(x) - 2 \langle Vh, u \Delta u\rangle + \sum_{x\in X}V  h^2 u^2(x)\\
& = \sum_{x\in X}V |\Delta u|^2(x) - 2 \langle Vh, |\nabla u|^2\rangle -\langle \Delta(Vh), u^2\rangle + \sum_{x\in X} V h^2 u^2(x)\\
& = \sum_{x\in X}V |\Delta u|^2(x) + \sum_{x\in X} \left[-\Delta(Vh) + Vh^2 + 2\frac{{\rm div}(Vh\nabla f)}{f}\right]u^2(x)\\
&\quad  - 2\F(Vh, f, u).
\end{align*}
Choose now $h = h_0 = \frac{\Delta f}{f} =: L(f)$ and using \eqref{Leib}, for any $x \in X$, there holds
\begin{align*}
& \quad f\left[-\Delta(Vh_0) + Vh_0^2 + 2\frac{{\rm div}(Vh_0\nabla f)}{f}\right](x) \\
& = \big[{-f\Delta(Vh_0) + Vh_0^2f - 2Vh_0\Delta f}\big](x) + 2 \langle\nabla f, \nabla(Vh_0)\rangle_x\\
& = -\big[f\Delta(Vh_0) + Vh_0\Delta f\big](x) + 2 \langle\nabla f, \nabla(Vh_0)\rangle_x\\
& = -\Delta(Vh_0f)(x) \\
& = -\Delta(V\Delta f)(x).
\end{align*}
Combining the above two equalities, we obtain a second order Hardy-Rellich equality:
\begin{align}\label{sec}
\sum_{x\in X}V |\Delta u|^2(x) = \sum_{x\in X} \frac{\Delta(V\Delta f)}{f}(x) u^2(x) + 2\F(VL(f), f, u) + \sum_{x\in X}V\left|\Delta u - L(f) u\right|^2(x).
\end{align}
Finally, using iteration, we arrive at
\begin{thm}
\label{thm2.1}
Let $m\geq 1$, $V, f\in C(X)$ satisfy $\Delta^k f > 0$ for all $0 \leq k \leq m-1$. Then for any $u\in C_c(X)$, there holds
\begin{align}
\label{equam}
\begin{split}
& \quad \sum_{x\in X}V |\Delta^m u|^2 - \sum_{x\in X} \frac{\Delta^m(V\Delta^m f)}{f} u^2\\
& = \sum_{k = 0}^{m-1}\sum_{x\in X}\frac{\Delta^{m-1-k}(V\Delta^m f)}{\Delta^{k+1} f}\left|\Delta^{k+1} u - L(\Delta^k f) \Delta^k u\right|^2(x)\\
& \quad + 2\sum_{k = 0}^{m-1}\F\left(\frac{\Delta^{m-1-k}(V\Delta^m f)}{\Delta^kf}, \Delta^k f, \Delta^k u\right).
\end{split}
\end{align}
\end{thm}

\medskip\noindent
{\bf Proof.} In fact, \eqref{sec} implies that \eqref{equam} holds for $m=1$. Assume \eqref{equam} for $m=j$, let we consider $m=j+1$. As $\Delta^j f > 0$ over $X$, using $\Delta^j u$ instead of $u$ and $\Delta^j f$ instead of $f$ in \eqref{sec}, we obtain
\begin{align*}
\sum_{x\in X}V |\Delta^{j+1} u|^2 & = \sum_{x\in X} \frac{\Delta(V\Delta^{j+1}f)}{\Delta^j f} |\Delta^j u|^2 + 2\F\left(VL(\Delta^j f), \Delta^j f, \Delta^j u\right)\\
& \quad + \sum_{x\in X}V\left|\Delta^{j+1} u - L(\Delta^j f) \Delta^j u\right|^2.
\end{align*}
Denote $W=\frac{\Delta(V\Delta^{j+1}f)}{\Delta^j f}$, applying \eqref{equam} with $W$ and $m =j$, there holds
\begin{align*}
\sum_{x\in X}V |\Delta^{j+1} u|^2 & =\sum_{x\in X} \frac{\Delta^j(W\Delta^j f)}{f} u^2+ \sum_{k = 0}^{j-1}\sum_{x\in X}\frac{\Delta^{j-1-k}(W\Delta^j f)}{\Delta^{k+1} f}\left|\Delta^{k+1} u - L(\Delta^k f) \Delta^k u\right|^2(x)\\
& \quad + 2\sum_{k = 0}^{j-1}\F\left(\frac{\Delta^{j-1-k}(W\Delta^j f)}{\Delta^kf}, \Delta^k f, \Delta^k u\right)+ 2\F\left(VL(\Delta^j f), \Delta^j f, \Delta^j u\right)\\
& \quad + \sum_{x\in X}V\left|\Delta^{j+1} u - L(\Delta^j f) \Delta^j u\right|^2\\
& = \sum_{x\in X} \frac{\Delta^{j+1}(V\Delta^{j+1} f)}{f} u^2 + \sum_{k = 0}^{j}\sum_{x\in X}\frac{\Delta^{j-k}(V\Delta^{j+1} f)}{\Delta^{k+1} f}\left|\Delta^{k+1} u - L(\Delta^k f) \Delta^k u\right|^2(x)\\
& \quad + 2\sum_{k = 0}^{j}\F\left(\frac{\Delta^{j-k}(V\Delta^{j+1} f)}{\Delta^kf}, \Delta^k f, \Delta^k u\right).
\end{align*}
So we are done. \qed

\medskip
If $\Delta^m f > 0$ over $X$, using \eqref{FH} with $\Delta^m u$, $\Delta^m f$ instead of $u$ and $f$ respectively, there holds
\begin{align*}
\sum_{x\in X} V(x)|\nabla \Delta^m u|^2(x) = \sum_{x\in X} \frac{\T(V, \Delta^m f)}{\Delta^m f}(x) |\Delta^m u|^2(x) + \F(V, \Delta^m f, \Delta^m u) .
\end{align*}
Applying Theorem \ref{thm2.1} to the first term on the RHS, we can claim immediately the following Hardy-Rellich equality of odd order.
\begin{thm}
\label{thm2.2}
Let $m \geq 1$, $V, f\in C(X)$ with $\Delta^k f > 0$ for $0\leq k\leq m$. Then for any $u\in C_c(X)$,
\begin{align*}
& \quad \sum_{x\in X} V|\nabla(\Delta^m u)|^2(x) - \sum_{x\in X}\frac{\Delta^m \T(V,\Delta^mf)}{f}u^2(x)\\
& = \sum_{k=0}^{m-1}\sum_{x\in X}\frac{\Delta^{m-1-k}\T(V,\Delta^mf)}{\Delta^{k+1}f}\left|\Delta^{k+1}u-L(\Delta^k f)\Delta^k u\right|^2(x) \\& \quad + 2 \sum_{k=0}^{m-1}\F\left(\frac{\Delta^{m-1-k}\T(V,\Delta^mf)}{\Delta^k f},\Delta^k f, \Delta^k u\right)+\F\left(V, \Delta^mf, \Delta^m u\right).
\end{align*}
\end{thm}

Here we keep the notation $\T$ to shorten the writing. Theorems \ref{thm2.1}-\ref{thm2.2} allow us to establish some Hardy-Rellich inequalities. Furthermore, the best constant or the best weights could be obtained with suitable choice of $f$ following the optimization on the sequel weights $\frac{\Delta^m(V\Delta^m f)}{f}$ or $\frac{\Delta^m\T(V, \Delta^m f)}{f}$.

\smallskip
In particular, let $V\equiv 1$, if $\Delta^k f>0$ for $0\leq k\leq \lfloor\frac{\ell-1}{2}\rfloor$ and $\Delta^i f\geq 0$ for $\lfloor\frac{\ell-1}{2}\rfloor<i\leq \ell-1$ with $\lfloor \beta\rfloor$ the maximum integer not exceeding $\beta$, there holds
\begin{align}
\label{rellich-f}
\sum_{x\in X}|\Delta^{\frac{\ell}{2}} u|^2(x) \geq \sum_{x\in X} \frac{\Delta^\ell f}{f} u^2(x), \quad \forall\; u\in C_c(X).
\end{align}

\begin{rem}
The idea of using equalities to handle Hardy-Rellich type inequalities was frequently used. The equality \eqref{FH} generalizes the ground state transform in \cite{KPP} (which is equivalent to the case of $V\equiv 1$) to the weighted case. The main toolbox in \cite[section 2]{KPP} contains the product and chain rules; the ground state transform and a discrete co-area formula, which are all based on some identities. Interested readers can also take a look at \cite{HY} for similar ideas in the continuous setting. For example, let $\alpha \in \R$ and $\varphi \in C_c^1(0, \infty)$, if $v = x^\frac{\alpha-1}{2}\varphi$, there holds
\begin{align*}
\int_0^\infty x^\alpha|\varphi'(x)|^2 dx & = \frac{(\alpha-1)^2}{4} \int_0^\infty x^{\alpha-2}\varphi^2(x) dx + \int_0^\infty x|v'(x)|^2 dx
\end{align*}
which gives the sharp weighted Hardy inequality over $(0, \infty)$ corresponding to \eqref{H} or \eqref{a<0}.
\end{rem}

\section{Hardy inequalities on $\N$ with shifting weights}
\reset
Here we will prove Theorem \ref{thm1.2}. As mentioned before, our departure point is the weighted equality \eqref{equa1}. Let $V_n\geq 0$, $f_n>0$, $\forall\; n\geq 1$ with $f_0=0$, we get immediately the first order Hardy inequality: For any $u \in C_c(\N)$ with $u_0 = 0$,
\begin{align}
\label{equa31}
\sum_{n=1}^\infty V_n|\nabla u_n|^2 \geq -\sum_{n=1}^\infty \frac{{\rm div}(V\nabla f)_n}{f_n} u^2_n, \quad \forall\; u \in C_c(\N), \; u_0 = 0.
\end{align}
The equality can not be attained if $V_n > 0$, since the RHS of \eqref{equa1} is zero, which implies that $u$ is parallel to $f$ over $\N$, this is impossible as $u \in C_c(\N)$. In the following, we will choose suitable $f_n$ to prove the first order estimates \eqref{a<0}.

\medskip
Seeing \eqref{a<0}, we take $V_n = (n-1)^\alpha$ if $n \geq 2$, $V_1 = 0$ and $f_n=n^{\frac{1-\alpha}{2}}$ on $\N$. Then for $n\geq 2$,
\begin{align*}
-\frac{{\rm div}(V \nabla f)_n}{f_n} & = \frac{V_n (f_n-f_{n-1}) - V_{n+1}(f_{n+1}-f_n)}{f_n}\\
& = V_n + V_{n+1} - \frac{f_{n-1}}{f_n}V_n - \frac{f_{n+1}}{f_n} V_{n+1}\\
& = n^\alpha\left[ 1+ \Big(1-\frac{1}{n}\Big)^\alpha- \Big(1-\frac{1}{n}\Big)^{\frac{1+\alpha}{2}}- \Big(1+\frac{1}{n}\Big)^{\frac{1-\alpha}{2}}\right] =: n^\alpha H_\alpha\left(\frac{1}{n}\right)
\end{align*}
where
$$
H_\alpha(x)=1+ (1-x)^\alpha-(1-x)^{\frac{1+\alpha}{2}}-(1+x)^{\frac{1-\alpha}{2}} \quad \mbox{on $[0,1)$.}
$$
Direct computation yields
\begin{align*}
\frac{H'''_\alpha(x)}{(1-\alpha)(1-x)^{\alpha-3}} & = \alpha(\alpha-2) + \frac{(1+\alpha)(3-\alpha)}{8}(1-x)^{\frac{1-\alpha}{2}}\\
& \quad -\frac{(1+\alpha)(3+\alpha)}{8}(1+x)^{-\frac{\alpha+5}{2}}(1-x)^{3-\alpha}\\
& =: I_1 +I_2 +I_3.
\end{align*}
As $\alpha<0$, so $I_1>0$. We claim that $H_\alpha'''\geq 0$ over $[0, 1)$ for any $\alpha < 0$.
\begin{itemize}
\item If $\alpha \leq -5$, $I_2 , I_3 \leq 0$ for $x \in [0, 1)$. Moreover, $$(1+x)^{-\frac{\alpha+5}{2}}(1-x)^{3-\alpha} = (1-x^2)^{-\frac{\alpha+5}{2}}(1-x)^{\frac{11-\alpha}{2}} \in (0, 1],$$
then
\begin{align}
\label{est3.1}
I_1 +I_2 + I_3 \geq I_1 + \frac{(1+\alpha)(3-\alpha)}{8} -\frac{(1+\alpha)(3+\alpha)}{8} = \frac{\alpha(3\alpha - 9)}{4}>0.
\end{align}
\item If $\alpha\in (-5, -3]$, we have always $I_2 , I_3 \leq 0$. As $-\frac{\alpha+5}{2} \leq 0$, \eqref{est3.1} is still valid.
\item Let $\alpha\in (-3,-1)$, then $I_2 < 0$ and $I_3 \geq 0$ for $x \in [0, 1)$, so that
$$
I_1+I_2 +I_3\geq -\alpha(2-\alpha) + \frac{(1+\alpha)(3-\alpha)}{8}=\frac{7(\alpha-1)^2 -4}{8} >0.
$$
\item When $\alpha\in [-1,0)$ and $x\in [0,1)$, there holds
\begin{align*}
I_2 +I_3 & = \frac{1+\alpha}{8}(1-x)^{\frac{1-\alpha}{2}}\left[(3-\alpha)-(3+\alpha)(1-x)^{\frac{5-\alpha}{2}}(1+x)^{-\frac{5+\alpha}{2}}\right]\\
& \geq \frac{1+\alpha}{8}(1-x)^{\frac{1-\alpha}{2}}(-2\alpha)\geq 0.
\end{align*}
\end{itemize}
Therefore, $H'''_\alpha(x) >0$ for any $\alpha < 0$ and any $x \in [0, 1)$. As $H'_\alpha(0)=0$, $H''_\alpha(0)=\frac{(1-\alpha)^2}{2}$, we get $H(x)>\frac{(1-\alpha)^2}{4} x^2$ in $(0, 1)$. Hence
$$
-\frac{{\rm div}(V \nabla f)_n}{f_n} = n^\alpha H_\alpha\left(\frac{1}{n}\right) \geq \frac{(\alpha-1)^2}{4}n^{\alpha-2} \quad \mbox{for all}\; n \geq 2.
$$
As $V_1 =0$ and $u_1 = 0$, using \eqref{equa31}, we obtain \eqref{a<0} for all $\alpha < 0$. Clearly, the inequality \eqref{a<0bis} can be derived by considering $u_n = (\tau_{-1}w)_n = w_{n-1}$ for $n \geq 1$.\qed

\medskip
Similarly, we can prove quickly the Leray type inequality \eqref{leray1}. Choose now $f_n= \sqrt{\ln n}$ with $n\geq 2$, $f_1 =\epsilon > 0$ and $f_0=0$.  For $n\geq 3$, let $h_n(t)=\sqrt{\ln (n+t)}$ on $[{-1}, 1]$ and
\begin{align*}
g_n(t) = 2h_n(0) - h_n(t) - h_n(-t), \quad \forall\; |t| \leq 1.
\end{align*}
Then $g$ is an even function, $g_n(0)=0$, $g'_n(0)=0$ and $g_n''(0)=-2h_n''(0)$. The Taylor expansion yields
\begin{align*}
\Delta f_n=g_n(1) & = -h''_n(0)+ \frac{1}{3!}\int_0^1(1-s)^3 g^{(4)}_n (s) ds \geq -h''_n(0) = \frac{1}{4 n^2(\ln n)^{\frac{3}{2}}} + \frac{1}{2 n^2(\ln n)^{\frac{1}{2}}}.
\end{align*}
The last estimate is given by $g_n^{(4)}(t) = -h_n^{(4)}(t) - h_n^{(4)}(-t)>0$. In fact, direct computation gives
\begin{align*}
\big(\sqrt{\ln s}\big)''' = \frac{1}{s^3\sqrt{\ln s}} + \frac{3}{4s^3(\ln s)^\frac{3}{2}} + \frac{3}{8s^3(\ln s)^\frac{5}{2}},
\end{align*}
which is readily decreasing in $(1, \infty)$ and $h_n^{(4)}(t) < 0$ for $|t| \leq 1$, if $n \geq 3$. As $V_n = n$,
\begin{align*}
-{\rm div}(V \nabla f)_n  = n \Delta f_n -\nabla f_{n+1} & \geq \frac{1}{4 n(\ln n)^{\frac{3}{2}}} + \frac{1}{2 n(\ln n)^{\frac{1}{2}}}-\frac{\ln(1+\frac{1}{n})}{\sqrt{\ln(n+1)}+\sqrt{\ln n}}\\
& \geq \frac{1}{4 n(\ln n)^{\frac{3}{2}}}.
\end{align*}
Here we used $\ln(1+x) < x$ for $x > 0$. So we get
\begin{align*}
-\frac{{\rm div}(V\nabla f)_n}{f_n}\geq \frac{1}{4 n(\ln n)^2}, \quad \forall\; n \geq 3.
\end{align*}
For $n=2$, there holds, for $\epsilon$ small enough,
\begin{align*}
-\frac{{\rm div}(V\nabla f)_2}{f_2}= V_2 + V_3 - \frac{f_1}{f_2}V_2 - \frac{f_3}{f_2} V_3 = 5- 3\sqrt{\frac{\ln 3}{\ln 2}} - \frac{2\epsilon}{\sqrt{\ln 2}} >\frac{1}{8(\ln2)^2}.
\end{align*}
As $u_0 = u_1 = 0$, the proof is completed seeing \eqref{equa1} or \eqref{equa31}. \qed

\section{Hardy-Rellich inequality of any order on $\N$}
\reset
Here we will apply the weighted first order Hardy inequalities to get higher order Rellich inequalities with optimal constants, i.e. Theorem \ref{thm1.1}. Recall that $\Delta = - \nabla({\rm div})$, let $v = {\rm div} u$, then $v_0 = 0$ if $u_0 = u_1 = 0$, there holds then
\begin{align*}
\sum_{n \geq 1} |\Delta u_n|^2  = \sum_{n \geq 1} |\nabla v_n|^2 & \geq \frac{1}{4}\sum_{n \geq 1} \frac{|v_n|^2}{n^2} = \frac{1}{4}\sum_{n \geq 2} \frac{|\nabla u_n|^2}{(n-1)^2} \geq \frac{9}{16}\sum_{n \geq 2}\frac{u_n^2}{n^4}.
\end{align*}
Here we used \eqref{a<0} with $\alpha =-2$ to get the last estimate. Hence \eqref{rellich-l} holds for $\ell = 2$.

\medskip
Now assume that \eqref{rellich-l} holds true for $\ell \leq 2m$ for some $m \geq 1$. We will prove \eqref{rellich-l} for $\ell = 2m+1$ and $2m+2$. Consider $u \in C_c(\N)$ with $u_0 = u_1 \ldots = u_{2m} = 0$. Denote $v = \tau_1u$, then $\nabla v_i = 0$ for $0 \leq i \leq 2m-1$. Using successively \eqref{rellich-l} with $\ell = 2m$ and \eqref{a<0} with $\alpha = -4m$, there holds
\begin{align*}
\sum_{n \geq m+1} |\nabla(\Delta^m u)_n|^2 = \sum_{n \geq m} |\nabla(\Delta^m v)_n|^2 & = \sum_{n \geq m} |(\Delta^m \nabla v)_n|^2\\
& \geq \Big[\frac{(4m)!}{4^{2m}(2m)!}\Big]^2\sum_{n \geq  2m} \frac{(\nabla v)_n^2}{n^{4m}}\\
& = \Big[\frac{(4m)!}{4^{2m}(2m)!}\Big]^2\sum_{n\geq 2m+1} \frac{(\nabla u)_n^2}{(n - 1)^{4m}}\\
& = \Big[\frac{(4m)!}{4^{2m}(2m)!}\Big]^2\sum_{n\geq 2} \frac{(\nabla u)_n^2}{(n - 1)^{4m}}\\
& \geq \Big[\frac{(4m)!}{4^{2m}(2m)!}\Big]^2\times \frac{(4m+1)^2}{4} \sum_{n \geq 2} \frac{u_n^2}{n^{4m+2}}\\
& = \Big[\frac{(4m+2)!}{4^{2m+1}(2m+1)!}\Big]^2 \sum_{n \geq 2m+1} \frac{u_n^2}{n^{4m+2}},
\end{align*}
which means that \eqref{rellich-l} is valid for $\ell = 2m+1$.

\medskip
Similarly, for $u$ satisfying $u_k = 0$ with $0 \leq k \leq 2m+1$, let $v = {\rm div}(u)$, then $v_i = 0$ for $i \leq 2m$. Applying the above estimate with $\ell = 2m+1$,
\begin{align*}
\sum_{n\geq m+1} |(\Delta^{m+1} u)_n|^2 & = \sum_{n\geq m+1} |\nabla(\Delta^m v)_n|^2\\
& \geq \Big[\frac{(4m+2)!}{4^{2m+1}(2m+1)!}\Big]^2 \sum_{n \geq 2m+1} \frac{v_n^2}{n^{4m+2}}\\
& = \Big[\frac{(4m+2)!}{4^{2m+1}(2m+1)!}\Big]^2 \sum_{n \geq 2} \frac{(\nabla u_n)^2}{(n - 1)^{4m+2}}\\
& \geq \Big[\frac{(4m+2)!}{4^{2m+1}(2m+1)!}\Big]^2 \frac{(4m+3)^2}{4} \sum_{n \geq 2m+2} \frac{u_n^2}{n^{4m+4}}\\
& = \Big[\frac{(4m+4)!}{4^{2m+2}(2m+2)!}\Big]^2 \sum_{n \geq 2m+2} \frac{u_n^2}{n^{4m+4}}.
\end{align*}
The second inequality comes from \eqref{a<0} with $\alpha = -4m-2$, and $u_k =0$ for $k \leq 2m+1$. We conclude then \eqref{rellich-l} for $\ell = 2m+2$. The proof of Theorem \ref{thm1.1} is completed. \qed

\medskip
Connecting the discrete Hardy inequality to the continuous setting, we claim the sharpness of \eqref{rellich-l}.
\begin{prop}
\label{sharp}
For any $\ell \geq 2$, the constant $\Big[\frac{(2\ell)!}{4^{\ell}\ell!}\Big]^2$ is optimal for \eqref{rellich-l}, equally we cannot replace it by any bigger constant such that the Hardy-Rellich inequality remains true for all $u\in C_c(\N)$ satisfying $u_k = 0$, $0 \leq k \leq \ell - 1$.
\end{prop}

\medskip\noindent
{\bf Proof.} The key observation is that the discrete Hardy-Rellich estimate will deduce the corresponding inequality for the continuous setting. More precisely, we shall show that if
\begin{align}
\label{sharp1}
\sum_{2n \geq \ell} |\Delta^\frac{\ell}{2}u_n|^2  \geq \lambda\sum_{n= \ell}^\infty \frac{u_n^2}{n^{2\ell}}, \quad \forall\; u \in C_c(\N), \; u_k = 0, 0 \leq k \leq \ell -1;
\end{align}
with a constant $\lambda > 0$, then there holds
\begin{align}
\label{sharp2}
\int_0^\infty |\varphi^{(\ell)}(x)|^2dx  \geq \lambda \int_0^\infty \frac{\varphi(x)^2}{x^{2\ell}}dx, \quad \forall\; \varphi \in C_c^\infty(0, \infty).
\end{align}
This connection is elementary, but as far as we are aware, we did not find an explicit statement of such flavor in the literature.

\medskip
Since the best constant to \eqref{sharp2} is $\Big[\frac{(2\ell)!}{4^{\ell}\ell!}\Big]^2$, we need only to make sure that \eqref{sharp1} implies \eqref{sharp2}.

\medskip
Consider just $\ell = 2$, other cases can be proved by the same way. Assume that \eqref{sharp1} is valid for $\ell = 2$ and some $\lambda > 0$. Take $\varphi \in C_c^\infty(0, 1)$. For any integer $M \geq 1$, consider $u_n = \varphi(n/M)$. Hence for $M$ large, $u_0 = u_1 =0$. By Taylor-Lagrange formula,
$$\Delta u_n = 2\varphi\Big(\frac{n}{M}\Big) - \varphi\Big(\frac{n-1}{M}\Big) - \varphi\Big(\frac{n+1}{M}\Big) = -\frac{1}{M^2} \varphi''\Big(\frac{n}{M}\Big) + r_n$$
where $|r_n|\leq \frac{\|\varphi'''\|_\infty}{3M^3}$. If we take $w_n = M^\frac{3}{2}u_n$, then as $M$ goes to infinity,
\begin{align*}
\sum_{n \geq 1}(\Delta w_n)^2 = M^3\sum_{n = 1}^M (\Delta u_n)^2 = \frac{1}{M} \sum_{n = 1}^M \varphi''\Big(\frac{n}{M}\Big)^2 + O\left(M^{-1}\right) \quad \longrightarrow \quad \int_0^1 |\varphi''(x)|^2 dx.
\end{align*}
On the other hand, as $M \to \infty$,
\begin{align*}
\sum_{n \geq 1} \frac{w_n^2}{n^4} = \frac{1}{M} \sum_{n = 1}^M \Big(\frac{n}{M}\Big)^{-4}\varphi\Big(\frac{n}{M}\Big)^2 \quad \longrightarrow \quad \int_0^1 \frac{\varphi(x)^2}{x^4} dx.
\end{align*}
Using \eqref{sharp1} on $w_n$, we deduce that \eqref{sharp2} holds for any $\varphi \in C_c^\infty(0, 1)$, hence for any $\varphi \in C_c^\infty(0, \infty)$ by scaling argument. So we are done.\qed

\begin{rem}
\label{opti2}
With similar consideration, we see that the constant $\frac{(\alpha -1)^2}{4}$ in \eqref{a<0} and \eqref{H} is also sharp.
\end{rem}

\begin{rem}
\label{zeros}
The conditions $u_k = 0$ for $k \leq \ell -1$ are fundamentally necessary to claim the Hardy-Rellich type inequality \eqref{rellich-l}. Fix again $\ell =2$, other cases can be handled similarly. Let $u_1 \ne 0$, seeing $\Delta = \nabla(\tau_1\nabla)$, the idea is to construct sequences $(u_n)$ such that the first order Hardy inequality fails for $v = \tau_1\nabla u$ as $v_0 = u_1 - u_0 \ne 0$. Let $M \geq 2$,
\begin{equation*}
w_n = \left\{ \begin{array}{ll}
0 & n = 0;\\
1+\frac{1}{2M} & \mbox{if }\; 1 \leq n \leq M;\\
1 -  \frac{n-M}{M} & \mbox{if }\; M+1 \leq n \leq 3M;\\
 -1 +  \frac{n-3M}{2M} & \mbox{if }\; 3M+1 \leq n \leq 5M;\\
0 &  \mbox{if }\; n > 5M.
\end{array}
\right.
\end{equation*}
We can check that $\sum_{1 \leq n \leq 5M} w_n = 0.$ Define $u_n = \sum_{1 \leq k\leq n} w_k$. Clearly, $w_n = (\nabla u)_n$ for $n \geq 1$, $u \in C_c(\N)$, $u_0 = 0$ and $u_1 = 1+ \frac{1}{2M}$. Hence $\Delta u_n = w_n - w_{n+1} = 0$ for $n > 5M$ or $1 \leq n < M$, and
$$\sum_{n \geq 1} |\Delta u_n|^2 = O\left(\frac{1}{M}\right), \quad \mbox{while }\; \sum_{n \geq 2} \frac{u_n^2}{n^4} \geq \sum_{n = 2}^M \frac{u_n^2}{n^4} \geq \sum_{n = 2}^M \frac{1}{n^2}.$$
Tending $M$ to infinity, we see that if $u_1 \ne 0$ is allowed, the Rellich inequality \eqref{rellich-l} cannot hold generally even if we replace the RHS constant by any $\epsilon > 0$.
\end{rem}

\begin{rem}
The above example means that in the sense of quadratic form on $\ell^2(\N^*)$, there is no $c > 0$ such that $\Delta^2 \geq V_c$ with $V_c(n) = \frac{c}{n^4}$. Here $\N^* = \{1,2,3\ldots\}$, and this answers a question in \cite{GKS}. Very recently, an optimal result is given in \cite[Theorem 1.1]{GKS2}, which yields in particular that there is no non-trivial weight $V \geq 0$ such that $\Delta^2 \geq V$ over $\ell^2(\N^*)$.
\end{rem}

\section{Sharp Hardy inequalities for weights $n^\alpha$ with $\alpha \geq 0$}
\reset
Now we will prove \eqref{H}. As already mentioned, Gupta proved \eqref{H} for $\alpha\in[0,1)\cup[5,\infty)$ and $\alpha \in 2\N$ with different methods. Here we give a unified approach for all $\alpha\in[0,+\infty)$. In fact, with $V_n = n^\alpha$ and suitable $f_n$, we will show that
\begin{align}
\label{equa51}-\frac{{\rm div}(V\nabla f)_n}{f_n} =  V_n + V_{n+1} - \frac{f_{n-1}}{f_n}V_n - \frac{f_{n+1}}{f_n} V_{n+1} \geq \frac{(\alpha - 1)^2}{4}n^{\alpha-2} \quad \mbox{for any }\; n \geq 1.
\end{align}
Then \eqref{equa1} or \eqref{equa31} implies the desired Hardy inequalities. We decompose the study into three subcases.

\subsection{For $\alpha \in [0, 1]$.} Take $f_n=n^{\frac{1-\alpha}{2}}$ with $n\geq 0$, then for $n \geq 1$,
\begin{align}
\label{equa52}
-\frac{{\rm div}(V \nabla f)_n}{f_n}= n^\alpha \left[1+\Big (1+\frac{1}{n}\Big)^\alpha- \Big(1+\frac{1}{n}\Big)^{\frac{\alpha+1}{2}}- \Big(1 - \frac{1}{n}\Big)^{\frac{1-\alpha}{2}}\right] =n^\alpha G_\alpha\left(\frac{1}{n}\right)
\end{align}
where
\begin{align*}
G_\alpha(x)=1+(1+x)^\alpha-(1+x)^{\frac{\alpha+1}{2}}-(1-x)^{\frac{1-\alpha}{2}},\quad\mbox{for }\; 0\leq x\leq 1.
\end{align*}
Direct calculation gives $G_\alpha(0)=0$, $G'_\alpha(0)=0$, $G''_\alpha(0)=\frac{(\alpha-1)^2}{2}$ and
\begin{align*}
G'''_\alpha(x) & = \alpha(\alpha-1)(\alpha-2)(1+x)^{\alpha-3} + \frac{(1 - \alpha^2)(\alpha-3)}{8}(1+x)^{\frac{\alpha-5}{2}}\\
& \quad + \frac{(1 - \alpha^2)(\alpha+3)}{8}(1-x)^{\frac{-5-\alpha}{2}}\\
& =:\sum_{k=1}^3\widetilde{J}_k.
\end{align*}
As $\alpha \in [0, 1]$ and $(1+x)^{\frac{\alpha-5}{2}} \leq 1 \leq (1-x)^{\frac{-5-\alpha}{2}}$ for $x\in [0,1]$, we get
\begin{align*}
\widetilde{J}_2 + \widetilde{J}_3 \geq \frac{\alpha(1-\alpha^2)}{4} \geq 0.
\end{align*}
Hence $G'''_\alpha(x) \geq 0$ in $[0, 1]$, which yields $G_\alpha(x)\geq \frac{(\alpha-1)^2}{4} x^2$ in $[0, 1]$. Seeing \eqref{equa52}, we conclude that \eqref{equa51} holds true.

\subsection{For $\alpha \geq 3$} In the following, we will choose the shifting sequence $f_n=(n+1)^\frac{1-\alpha}{2}$ for $n\geq 1$ and $f_0=0$.

\medskip
First, $2^\alpha > (\alpha - 1)^2$ for $\alpha\geq 3$, we have then
\begin{align}
\label{est51}
-\frac{{\rm div}(V \nabla f)_1}{f_1} = V_1 + \Big(1 - \frac{f_2}{f_1}\Big)V_2 & = 1 + \Big[1 -\Big (\frac{2}{3}\Big)^\frac{\alpha -1}{2} \Big]2^\alpha \geq 1 + \frac{2^\alpha}{3} > \frac{(\alpha - 1)^2}{4}.
\end{align}

Let $n \geq 2$,
\begin{align}
\label{equa53}
\begin{split}
-\frac{{\rm div}(V \nabla f)_n}{f_n} & = n^\alpha \left[\Big(1+\frac{1}{n}\Big)^\alpha- \Big(1+\frac{1}{n}\Big)^{\frac{3\alpha-1}{2}}\Big(1+\frac{2}{n}\Big)^\frac{1-\alpha}{2} +1- \Big(1+\frac{1}{n}\Big)^\frac{\alpha-1}{2}\right]\\
& =: n^\alpha F_\alpha\Big(1 + \frac{1}{n}\Big)
\end{split}
\end{align}
where
\begin{align}
\label{Falpha}
F_\alpha(x) = 1+ x^\alpha -x^{\frac{3\alpha-1}{2}} (2x-1)^{\frac{1-\alpha}{2}}-x^{\frac{\alpha-1}{2}} = 1+ x^\alpha -x^{\frac{\alpha-1}{2}}-x^\alpha Y(x)^{\frac{\alpha-1}{2}}.
\end{align}
So we are interested in $F_\alpha$ over $I_0 := \big[1, \frac{3}{2}\big]$ with $Y(x)=\frac{x}{2x-1}\in[\frac{3}{4}, 1]$. Using $Y'(x)=-x^{-2}Y^2$, we can check that
\begin{align*}
F'''_\alpha(x)& = \alpha(\alpha-1)(\alpha-2)x^{\alpha-3}\Big[1-Y^{\frac{\alpha-1}{2}}\Big]-\frac{(\alpha-1)(\alpha-3)(\alpha-5)}{8}x^{\frac{\alpha-7}{2}}\\&
\quad +\frac{3(\alpha-2)(\alpha-1)}{4}x^{\alpha-4}Y^{\frac{\alpha+1}{2}}\Big[2(\alpha-1)-(\alpha+1)x^{-1}Y\Big]+\frac{(\alpha^2-1)(\alpha+3)}{8}x^{\alpha-6}Y^{\frac{\alpha+5}{2}}\\
& =:\sum_{k=1}^4 \widehat J_k.
\end{align*}
Moreover, $x^{-1}Y =\frac{1}{2x-1} \in \big[\frac{1}{2}, 1\big]$ and $xY = \frac{x^2}{2x -1} \geq 1$ over $I_0$, there holds then $$2(\alpha-1)-(\alpha+1)x^{-1}Y \geq \alpha -3,$$ so that
\begin{align*}
\widehat J_3 \geq \frac{3(\alpha-1)(\alpha-2)(\alpha-3)}{4}x^{\alpha-4}Y^{\frac{\alpha+1}{2}} & = \frac{3(\alpha-1)(\alpha-2)(\alpha-3)}{4}x^{\frac{\alpha-7}{2}}(xY)^{\frac{\alpha-1}{2}}Y
\\ & \geq \frac{9(\alpha-1)(\alpha-2)(\alpha-3)}{16}x^{\frac{\alpha-7}{2}}.
\end{align*}
Therefore
\begin{align*}
(\widehat J_3+ \widehat J_2 )x^{\frac{7-\alpha}{2}} \geq \frac{(\alpha-1)(\alpha-3)(7\alpha-8)}{16}\geq 0.
\end{align*}
It means that $F'''_\alpha(x)\geq 0$ in $I_0$ for $\alpha \geq 3$, since $\widehat J_1, \widehat J_4 \geq 0$ over $I_0$.

\medskip
With $F_\alpha(1)=0$,  $F'_\alpha(1)=0$, $F''_\alpha(1) =\frac{(\alpha-1)^2}{2}$, we get $F_\alpha(x) \geq \frac{(\alpha-1)^2}{4}(x-1)^2$ in $I_0$ by Taylor expansion. Combining with \eqref{est51} and \eqref{equa53}, we get \eqref{equa51} for $\alpha \geq 3$.

\subsection{For $\alpha \in (1,3)$} As above, we choose $f_n = (n+1)^{\frac{1-\alpha}{2}}$ for $n \geq 1$ and $f_0=0$. For $\alpha \in (1, 3)$,
\begin{align}
\label{est51bis}
-\frac{{\rm div}(V \nabla f)_1}{f_1} = 1 + \Big[1 -\Big (\frac{2}{3}\Big)^\frac{\alpha -1}{2} \Big]2^\alpha > 1 > \frac{(\alpha - 1)^2}{4}.
\end{align}
Let $F_\alpha$ be given by \eqref{Falpha}, $Y(x)$, $I_0$ be also as in the previous case. Notice that $x^{-1}Y=2Y-1$, we rewrite
\begin{align*}
F'''_\alpha(x) & = (\alpha-1)x^{\alpha-3}\left[\alpha(\alpha-2)-\frac{(\alpha-3)(\alpha-5)}{8}x^{\frac{-1-\alpha}{2}}+ L_\alpha(Y) \right]\\
& =: (\alpha-1)x^{\alpha-3}K_\alpha(x)
\end{align*}
with $Y^{\frac{1-\alpha}{2}}L_\alpha(Y) =$
\begin{align*}
(\alpha+1)(\alpha+3) Y^3 -\frac{3(3\alpha-1)(\alpha+1)}{2}Y^2 + \frac{9(3\alpha-1)(\alpha-1)}{4} Y - \frac{3(3\alpha-1)(3\alpha-5)}{8}.
\end{align*}
Consider $L_\alpha$ as function of variable $Y$, direct computation gives $L'_\alpha(Y)= Y^{\frac{\alpha-3}{2}}J_\alpha(Y)$,
where
\begin{align*}
J_\alpha(Y) & = \frac{(\alpha+1)(\alpha+3)(\alpha+5)}{2}Y^3 - \frac{3(3\alpha-1)(\alpha+1)(\alpha+3)}{4}Y^2 + \frac{9(3\alpha-1)(\alpha^2-1)}{8}Y\\
& \quad  -\frac{3(\alpha-1)(3\alpha-1)(3\alpha-5)}{16}.
\end{align*}
For $Y \in I_1 :=[\frac{3}{4}, 1]$ and $\alpha \in (1, 3)$, there holds
\begin{align*}
J'_\alpha(Y) & = \frac{3(\alpha+1)}{8}Y\left[4(\alpha+3)(\alpha+5)Y - 4(3\alpha-1)(\alpha+3) + 3(3\alpha-1)(\alpha-1)Y^{-1}\right]\\
& \geq \frac{3(\alpha+1)}{8}Y\Big[3(\alpha+3)(\alpha+5)  -4(3\alpha-1)(\alpha+3)+3(3\alpha-1)(\alpha-1)\Big]\\
& = \frac{3(\alpha+1)Y}{8}\times 20(3 - \alpha) > 0.
\end{align*}
Meanwhile, we have
\begin{align*}
J_\alpha\Big(\frac{3}{4}\Big) = \frac{3(-9\alpha^3+63\alpha^2-175\alpha+265)}{128} =: g(\alpha)
\end{align*}
Elementary calculus shows that $g$ is decreasing and positive in $[1, 3]$, Hence $J_\alpha(Y)>0$ for any $Y \in I_1$ and $\alpha\in (1, 3)$, so $L_\alpha$ is increasing in $I_1$. Moreover, as $\alpha\in (1, 3)$, using $xY \geq 1$, $Y\geq \frac{3}{4}$ and $J_\alpha(Y)\geq J_\alpha\left(\frac{3}{4}\right)$ for $x \in I_0$, there holds
\begin{align*}
K'_\alpha(x)& = \frac{(\alpha+1)(\alpha-3)(\alpha-5)}{16} x^{-\frac{3+\alpha}{2}} + Y^{\frac{\alpha-3}{2}}J_\alpha(Y) Y'(x)\\
& = x^{-\frac{3+\alpha}{2}}\left[\frac{(\alpha+1)(\alpha-3)(\alpha-5)}{16} - (xY)^{\frac{\alpha-1}{2}} Y J_\alpha(Y)\right]\\
& \leq x^{-\frac{3+\alpha}{2}} \left[\frac{(\alpha+1)(\alpha-3)(\alpha-5)}{16}-\frac{3}{4} J_\alpha\left(\frac{3}{4}\right)\right]\\
& =: x^{-\frac{3+\alpha}{2}} G(\alpha),
\end{align*}
where
\begin{align*}
G(\alpha) =\frac{(\alpha+1)(\alpha-3)(\alpha-5)}{16}-\frac{3}{4} g(\alpha)
 = \frac{113\alpha^3 -791 \alpha^2 +1799 \alpha -1905}{512}.
\end{align*}
We can check that $G(\alpha) <0$ for $\alpha \in (1,3)$, which implies then $K'_\alpha(x)<0$ in $I_0$. Moreover,
$$
K_\alpha(1) = \frac{3(\alpha^2-3\alpha+4)}{4}>0.
$$
Recall that $F'''_\alpha(x) = (\alpha -1)x^{\alpha-3}K_\alpha(x)$, so $F'''_\alpha$ changes eventually at most once its sign over $I_0$ from positive value to negative value, hence
\begin{align}\label{1}
\min_{I_0} F_\alpha''(x)= \min\Big\{ F''_\alpha(1),~F''_\alpha\left(\frac{3}{2}\right)\Big\}.
\end{align}
Clearly $F''_\alpha(1) =\frac{(\alpha-1)^2}{2} > 0$ and
\begin{align}\label{2}
F''_\alpha\left(\frac{3}{2}\right)= 
(\alpha-1)\left[\alpha\left(\frac{3}{2}\right)^{\alpha-2}-(\alpha+1)\left(\frac{3}{2}\right)^{\frac{3\alpha-1}{2}}2^{\frac{-3-\alpha}{2}}
-\frac{\alpha-3}{4}\left(\frac{3}{2}\right)^{\frac{\alpha-5}{2}}\right] =: Q(\alpha)
\end{align}
 We will check that
\begin{align}\label{3}
Q(\alpha) > \frac{(\alpha-1)^2}{2}\quad \mbox{for }~\alpha\in(1,3).
\end{align}
Indeed, we only need to prove $H(\alpha)>0$ in $(1, 3)$, where
$$
H(\alpha):=\alpha-(\alpha+1)\left(\frac{3}{4}\right)^{\frac{3+\alpha}{2}}-\frac{\alpha-3}{4}\left(\frac{2}{3}\right)^{\frac{\alpha+1}{2}}-\frac{\alpha-1}{2}\left(\frac{3}{2}\right)^{2-\alpha}.
$$
Directly calculus gives
\begin{align*}
H'(\alpha)& = 1-\left(\frac{3}{4}\right)^{\frac{3+\alpha}{2}}-\frac{(\alpha+1)}{2}\left(\frac{3}{4}\right)^{\frac{3+\alpha}{2}}\ln\frac{3}{4}
-\frac{1}{4}\left(\frac{2}{3}\right)^{\frac{\alpha+1}{2}}-\frac{\alpha-3}{8}\left(\frac{2}{3}\right)^{\frac{\alpha+1}{2}}\ln\frac{2}{3}\\
&\quad -\frac{1}{2}\left(\frac{3}{2}\right)^{2-\alpha}+\frac{\alpha-1}{2}\left(\frac{3}{2}\right)^{2-\alpha}\ln\frac{3}{2}
\end{align*}
and
\begin{align*}
H''(\alpha)=&-\left(\frac{3}{4}\right)^{\frac{\alpha+3}{2}}\ln\frac{3}{4}-\frac{\alpha+1}{4}(\ln\frac{3}{4})^2\left(\frac{3}{4}\right)^{\frac{\alpha+3}{2}}+\frac{1}{4}\ln\frac{3}{2}\left(\frac{2}{3}\right)^{\frac{\alpha+1}{2}}\\
&-\frac{\alpha-3}{16}(\ln\frac{3}{2})^2\left(\frac{2}{3}\right)^{\frac{\alpha+1}{2}}+\ln\frac{3}{2}\left(\frac{3}{2}\right)^{2-\alpha}-\frac{\alpha-1}{2}(\ln\frac{3}{2})^2\left(\frac{3}{2}\right)^{2-\alpha}\\
>&-\ln\frac{3}{4}\left(\frac{3}{4}\right)^{\frac{\alpha+3}{2}}\left(1+\ln\frac{3}{4}\right) + \ln\frac{3}{2}\left(\frac{3}{2}\right)^{2-\alpha}\left(1-\ln\frac{3}{2}\right)>0.
\end{align*}
Hence $H'$ is increasing. As $H'(1.5)\approx-0.0198<0$, $H'(1.55)\approx 0.0117>0$, there exists $1.5<\alpha_0< 1.55$ such that $H'(\alpha_0)=0$. The convexity of $H$ implies
$$
H(\alpha)\geq H(\alpha_0)>1.5-2.55\left(\frac{3}{4}\right)^{\frac{9}{4}}+\frac{1.45}{4}\left(\frac{2}{3}\right)^{\frac{2.55}{2}}-\frac{0.55}{2}\left(\frac{3}{2}\right)^{0.5}\approx 0.0445>0.
$$
Combining with \eqref{1}-\eqref{3}, we see that $F''_\alpha(x)\geq \frac{(\alpha-1)^2}{2}$ in $I_0$ for $\alpha\in (1,3)$. Again, by Taylor expansion at $x_0=1$, we get $F_\alpha(x) \geq \frac{(\alpha-1)^2}{4}(x-1)^2$. Using \eqref{est51bis}, \eqref{equa51} is valid for $\alpha \in (1, 3)$.

\medskip
To conclude, for any $\alpha \geq 0$, the weighted Hardy inequality \eqref{H} holds true. \qed

\medskip
\section{First order weighted Hardy inequality over $\mathbb{Z}^d$}
\reset
Here we attempt to apply our idea to consider the lattice $\Z^d$ with $d \geq 2$. By definition, $x\sim y$ on $\mathbb{Z}^d$ if and only if $y=x\pm e_i$ for $1\leq i \leq d$, and $b(x, y) = 1$ if $x \sim y$, here $\{e_i\}$ means the canonical basis in $\Z^d$. We claim the following weighted Hardy inequality, where $|\cdot|$ stands for the Euclidean norm in $\R^d$.

\begin{thm}
Let $\alpha>2-d$, there exists $A_{d,\alpha}\in\mathbb{R}$ such that
\begin{align}
\label{HZd}
\sum_{x\in \mathbb{Z}^d}|x|^\alpha|\nabla u|^2(x)\geq \sum_{x\in \mathbb{Z}^d} w(x) u^2(x),\quad \forall ~u\in C_c(\mathbb{Z}^d\setminus{\{0}\}),
\end{align}
where as $|x|\to +\infty$
$$w(x) = \frac{(d-2+\alpha)^2}{4}|x|^{\alpha-2} + A_{d,\alpha}|x|^{\alpha-4} + O(|x|^{\alpha -6}).$$
\end{thm}

\medskip\noindent
{\bf Proof.}  Take $f(x)=|x|^{2\gamma}$ and $V=|x|^\alpha$ for $x\in \mathbb{Z}^d\setminus{\{0}\}$, by the first order Hardy equality \eqref{FHdiv}, the estimate \eqref{HZd} holds true with $w(x) := \frac{-{\rm div}(V\nabla f)}{f}$ and $w(0)=0$. We will estimate $w(x)$ then make a suitable choice of $\gamma$. To avoid confusion, we use $\Delta f_x$ to denote the discrete Laplacian of $f$ at $x \in \Z^d$,
\begin{align*}
\frac{\Delta f_x}{f(x)} = \frac{1}{f(x)}\sum_{y\sim x}(f(x)-f(y)) & = \sum_{i=1}^d \left[2- \Big(1+\frac{2x_i}{|x|^2}+\frac{1}{|x|^2}\Big)^{\gamma}- \Big(1-\frac{2x_i}{|x|^2}+\frac{1}{|x|^2}\Big)^{\gamma}\right].
\end{align*}
By Taylor expansion, as $|x|\to +\infty$, there holds
\begin{align*}
\Big(1+\frac{1+2x_i}{|x|^2}\Big)^\gamma+ \Big(1+\frac{1-2x_i}{|x|^2}\Big)^\gamma & = 2+\gamma\frac{2}{|x|^2}+\frac{\gamma(\gamma-1)}{2}\frac{2+8x_i^2}{|x|^4}
+\frac{\gamma(\gamma-1)(\gamma-2)}{3}\frac{1+12x_i^2}{|x|^6} \\& \quad + \frac{\gamma(\gamma-1)(\gamma-2)(\gamma-3)}{12}\frac{1+24x_i^2+16x_i^4}{|x|^8}+O(|x|^{-6}),
\end{align*}
which yields
\begin{align*}
\frac{\Delta f_x}{f(x)} & = \sum_{i=1}^d \left[2-(1+\frac{2x_i}{|x|^2}+\frac{1}{|x|^2})^{\gamma}-(1-\frac{2x_i}{|x|^2}+\frac{1}{|x|^2})^{\gamma}\right]\\
& = \frac{-2d\gamma-4\gamma(\gamma-1)}{|x|^2} - \frac{\gamma(\gamma-1)[d+4(\gamma-2)]}{|x|^4}-\frac{4\gamma(\gamma-1)(\gamma-2)(\gamma-3)}{3}\sum_{i=1}^d \frac{x_i^4}{|x|^8}\\&\quad +O(|x|^{-6}).
\end{align*}
Meanwhile, we have
\begin{align*}
\frac{\langle \nabla f, \nabla V\rangle_x}{f(x)}& = \frac{1}{2f(x)}\sum_{x\sim y} \big[f(x)-f(y)\big]\big[V(x)-V(y)\big]\\
& = \frac{1}{2f(x)}\sum_{i=1}^d \Big[|x|^{2\gamma}-(|x|^2+2x_i+1)^\gamma\Big]\Big[|x|^\alpha-(|x|^2+2x_i+1)^{\frac{\alpha}{2}}\Big]\\
&\quad +\frac{1}{2f(x)}\sum_{i=1}^d \Big[|x|^{2\gamma}-(|x|^2-2x_i+1)^\gamma\Big]\Big[|x|^\alpha-(|x|^2-2x_i+1)^{\frac{\alpha}{2}}\Big]\\
& = \frac{|x|^\alpha}{2}\sum_{i=1}^d (J_+ + J_-),\\
\end{align*}
where
\begin{align*}
J_\pm: = \Big[1- \Big(1 \pm 2\frac{x_i}{|x|^2}+\frac{1}{|x|^2}\Big)^\gamma\Big] \Big[1- \Big(1 \pm 2\frac{x_i}{|x|^2}+\frac{1}{|x|^2}\Big)^{\frac{\alpha}{2}}\Big].
\end{align*}
As $|x|\to +\infty$, direct expansions yield
\begin{align*}
J_\pm = \frac{\alpha\gamma}{2}\frac{(1 \pm 2x_i)^2}{|x|^4}+\frac{\alpha\gamma(\alpha+2\gamma-4)}{8}\frac{(1\pm 2x_i)^3}{|x|^6} + K_{\alpha,\gamma}\frac{(1\pm 2x_i)^4}{|x|^8}+O(|x|^{-6})
\end{align*}
with
\begin{align*}
K_{\alpha,\gamma}= \frac{\alpha\gamma}{48}\Big[(\alpha-2)(\alpha-4)+3(\alpha-2)(\gamma-1)+4(\gamma-1)(\gamma-2)\Big].
\end{align*}
Therefore, as $|x|$ goes to infinity, there holds
\begin{align*}
\frac{1}{2}\sum_{i=1}^d (J_+ + J_-) = \frac{2\alpha\gamma}{|x|^2}+ \frac{\alpha\gamma(d+3(\alpha+2\gamma-4))}{2|x|^4} + 32K_{\alpha,\gamma}\sum_{i=1}^d \frac{x_i^4}{|x|^8}+O(|x|^{-6}).
\end{align*}
This implies as $|x|\to+\infty$,
\begin{align*}
\frac{\omega(x)}{V(x)}& = \frac{\Delta f_x}{f(x)} - \frac{\langle \nabla f, \nabla V\rangle_x}{f(x)V(x)}\\
& =  \frac{-2\gamma(d + \alpha - 2 +2\gamma)}{|x|^2}-\gamma\left[(\gamma-1)(d+4(\gamma-2))+\frac{d\alpha}{2}+\frac{3\alpha(\alpha+2\gamma-4)}{2}\right]\frac{1}{|x|^4}\\&
\quad - \left[\frac{4\gamma(\gamma-1)(\gamma-2)(\gamma-3)}{3}+32K_{\alpha,\gamma}\right]\sum_{i=1}^d \frac{x_i^4}{|x|^8}+O(|x|^{-6}).
\end{align*}
Now we choose $\gamma=\frac{2-d-\alpha}{4} < 0$ to optimize the coefficient of the main order term, then
\begin{align*}
\omega(x) & = \frac{(d-2+\alpha)^2}{4} |x|^{\alpha-2} + \frac{(d-2+\alpha)(3d+6+2\alpha^2-5\alpha)}{8}|x|^{\alpha - 4}\\
& \quad  - C_{\alpha, d}
\sum_{i=1}^d\frac{x_i^4}{|x|^{8-\alpha}}+O(|x|^{\alpha-6})\\
& \geq \frac{(d-2+\alpha)^2}{4}\frac{1}{|x|^{2-\alpha}}+ A_{\alpha,d}\frac{1}{|x|^{4-\alpha}}+O(|x|^{\alpha-6}),
\end{align*}
since $C_{\alpha, d} \sum x_i^4 \leq \max(C_{\alpha, d}, 0)|x|^4$.\qed

\begin{rem}
\label{rem5.1}
Let $d \geq 3$, $\alpha=0$, i.e $V\equiv 1$, it was shown in \cite[Theorem 7.2]{KPP} that \eqref{HZd} holds with
$$w(x) = \frac{(d-2)^2}{4}\frac{1}{|x|^2} + O(|x|^{-3}),  \quad \mbox{as }\; |x| \to \infty.$$
They used the asymptotic behavior of the Green's function for discrete Laplacian over $\Z^d$. We don't know if their technique works for weighted situation with $\alpha \ne 0$, our approach seems to be more computable and provides a weight $w$ with more precise behavior. In fact, for $\alpha=0$, we have
\begin{align*}
\omega(x) & = \frac{(d-2)^2}{4}\frac{1}{|x|^2}+\frac{3(d^2-4)}{8}\frac{1}{|x|^4}-\frac{(d^2-4)(d+6)(d+10)}{192}\sum_{i=1}^d\frac{x_i^4}{|x|^8}+O(|x|^{-6})\\
& \geq \frac{(d-2)^2}{4}\frac{1}{|x|^2}-\frac{(d^2-4)(d^2+16d-12)}{192}\frac{1}{|x|^4}+O(|x|^{-6}).
\end{align*}
\end{rem}

As for Theorem \ref{leray}, we can also get Leray type estimate in $\Z^d$. For simplicity, we consider just $\Z^2$ case here, that is
\begin{align}
\label{lerayZ2}
\sum_{x\in\mathbb{Z}^2}|\nabla u|^2(x)\geq \sum_{x\in\mathbb{Z}^2} w_2(x) u^2(x), \quad \forall \; u \in  C_c(\Z^2\cap{\{|x|\geq 2}\})
\end{align}
where
$$
w_2(x) = \frac{1}{4|x|^2(\ln|x|)^2} + \frac{2(x_1^4+x_2^4)}{|x|^4}- \frac{3}{2|x|^4\ln|x|} +O\left( \frac{1}{|x|^4\ln|x|^2}\right)\quad \mbox{as } |x|\to + \infty.
$$
Let $f=(\ln|x|)^\gamma$ with $|x|\geq 2$, and consider
$$h_i(t)=(\ln|x+te_i|)^\gamma, \;\; i=1,2 \quad \mbox{and}\quad  g_i(t) =h_i(t)+h_i(-t)-2h_i(0).$$ We have then
$$
g_i(0)= g'_i(0)= g'''_i(0)=0,\quad g''_i(0)=2h''_i(0)=2\frac{\partial^2 f}{\partial x_i^2}(x),
$$
and $\Delta f_x= - g_1(1) - g_2(1)$. As before, Taylor expansion yields
\begin{align*}
\Delta f_x & =-\sum_{i=1}^2\frac{\partial^2 f}{\partial x_i^2}(x) - \frac{g_1^{(4)}(0) + g_2^{(4)}(0)}{24} - \frac{1}{120}\sum_{i=1}^2\int_0^1(1-s)^5 g^{(6)}_i(s)ds\\
& = \frac{-\gamma(\gamma-1)}{|x|^2(\ln|x|)^{2-\gamma}} - \frac{g_1^{(4)}(0) + g_2^{(4)}(0)}{24} + J.
\end{align*}
Moreover, it's not difficult to see that $|h^{(6)}_i(t)| = O(|x|^{-6}(\ln|x|)^{\gamma-1})$ as $|x|\to+\infty$ uniformly for $t\in[-1,1]$, which means
$J = O(|x|^{-6}(\ln|x|)^{\gamma-1})$ as $|x|\to+\infty$.

\medskip
Set now $\gamma=\frac{1}{2}$, direct computation yields
\begin{align*}
\frac{1}{f(x)}\sum_{i=1}^2g^{(4)}_i(0) & = \left[36-48\sum_{i=1}^2\frac{x_i^4}{|x|^4}\right]\frac{1}{|x|^4\ln|x|} -\frac{1}{2}\left[26\sum_{i=1}^2\frac{x_i^4}{|x|^4}-21\right]\frac{1}{|x|^4(\ln|x|)^2}\\
&\quad +\frac{3}{4}\left[6-9\sum_{i=1}^2\frac{x_i^4}{|x|^4}\right]\frac{1}{|x|^4(\ln|x|)^3}-\frac{15}{8}\frac{1}{|x|^4(\ln|x|)^4}\sum_{i=1}^2\frac{x_i^4}{|x|^4}\\
& = \left[36-48\sum_{i=1}^2\frac{x_i^4}{|x|^4}\right]\frac{1}{|x|^4\ln|x|} +O(|x|^{-4}(\ln|x|)^{-2}).
\end{align*}
Applying \eqref{HZd} with $V \equiv 1$, we get a Hardy inequality \eqref{lerayZ2} with
\begin{align*}
w_2(x) = \frac{\Delta f_x}{f(x) } & = \frac{1}{4}\frac{1}{|x|^2(\ln|x|)^2}+ \frac{1}{24}\left[48\sum_{i=1}^2\frac{x_i^4}{|x|^4}-36\right]\frac{1}{|x|^4\ln|x|} +O(|x|^{-4}(\ln|x|)^{-2})\\
& \geq \frac{1}{4}\frac{1}{|x|^2(\ln|x|)^2}-\frac{1}{2} \frac{1}{|x|^4\ln|x|} +O(|x|^{-4}(\ln|x|)^{-2}).
\end{align*}

\section{Further remarks}
\reset
Using the weighted equality \eqref{equa1} or more generally \eqref{FH}, we establish some Hardy-Rellich inequalities with optimal constants in Theorems \ref{thm1.1} to \ref{prop1.2'}. Our approach can be adapted to the $\ell^p$ setting. There exist however many rooms for improvement for higher order Rellich inequalities over $\N$, and more mysteries exist for the lattice situation $\Z^d$ with $d\geq 2$, even for the first order Hardy inequality.

\subsection{} Firstly, we want to point out that the idea leading to the Hardy inequality \eqref{FHdiv} on graphs works also for general $\ell^p$ setting. Let $(X, b)$ be a locally finite graph. For $u \in C(X)$, denote
\begin{align}
|\nabla u|^p_p(x) := \frac{1}{2} \sum_{y\sim x} b(x,y)^{p-1}|\nabla u(x, y)|^p, \quad \forall \; x \in X.
\end{align}
The above definition can be motivated as follows. Let $b(x,y)^{-1}$ mean the direct distance between vertices $x, y$, then the unitary gradient of $u$ is given by $b(x,y)\nabla u(x, y)$, therefore the accumulated pointwise $\ell^p$ norm at $x$ is given by $|\nabla u|^p_p(x)$.

\medskip
For simplicity, we use the following notation
$$t^\beta = |t|^{\beta - 1}t, \quad \forall \; \beta > 0, \; t \in \R.$$
We define also
\begin{align}
(\nabla u)^{p-1}_b(x, y) := b(x,y)^{p-2}\nabla u(x, y)^{p-1}{\mathbbm 1}_{x\sim y}, \quad \forall \;x, y \in X, u \in C(X).
\end{align}
Remark that $(\nabla u)^p$ is anti-symmetric over $X\times X$ and ${\rm div}\big[(\nabla u)^{p-1}_b\big] = \Delta_p u$.

\medskip
Let $f\in C(X)$ be a positive function. By the discrete Picone's inequality (see \cite{A}), there holds
\begin{align*}
|\nabla u(x, y)|^p \geq \nabla\left(\frac{u^p}{f^{p-1}}\right)(x,y) (\nabla f)^{p-1}(x,y), \quad \forall\; x, y \in X, u\in C(X).
\end{align*}
Therefore, let $u \in C_c(X)$,
\begin{align*}
\sum_{x\in X} V(x) |\nabla u|^p_p(x) & = \frac{1}{2}\sum_{x,y\in X}V(x) b(x,y)^{p-1}|\nabla u(x, y)|^p
\\ & \geq \frac{1}{2} \sum_{x,y\in X} b(x,y)^{p-1} V(x) \nabla\left(\frac{u^p}{f^{p-1}}\right)(x,y)\nabla f(x, y)^{p-1}\\
& = \frac{1}{2}\sum_{x,y\in X} b(x,y)^{p-1} \big[V(x) + V(y)\big] \frac{u^p}{f^{p-1}}(x)\nabla f(x, y)^{p-1}\\
& = - \sum_{x\in X}\frac{{\rm div}\big[V (\nabla f)_b^{p-1}\big]}{f^{p-1}}(x) u^{p-1}(x).
\end{align*}
For the last two lines above, we used the symmetry of $b$ and the anti-symmetry of $(\nabla f)^{p-1}$. The divergence operator is that given in \eqref{div}. Then for any $p > 1$, we get the $\ell^p$- discrete Hardy inequality over $X$, which generalizes clearly \eqref{FHdiv} with $p = 2$.
\begin{align}
\label{hardyp}
\sum_{x\in X} V(x) |\nabla u|^p_p(x)\geq - \sum_{x\in X}\frac{{\rm div}\big[V (\nabla f)_b^{p-1}\big]}{f^{p-1}}(x) u^{p}(x), \quad \forall\; u \in C_c(X).
\end{align}
The proof of \cite[Proposition 3]{FKP} used somehow an equivalent form of \eqref{hardyp} with $V \equiv 1$ over $\N$, see also \cite[Theorem 3.1]{Fi} for locally summable graphs.

\subsection{} As Keller-Pinchover-Pogorzelski's result in \eqref{kpp1} for $\ell = 1$, we can obtain sharper weights by adding infinitely many positive remainder terms in \eqref{rellich-l}. Consider for example $\ell = 2$: Let $u \in C_c(\N)$ with $u_0 = u_1 = 0$,
seeing the proof of \eqref{a<0} for $\alpha = -2$, we arrive at
\begin{align}
\label{a=-2bis}
\sum_{n \geq 2} (n-1)^{-2} |\nabla u_n|^2  \geq \sum_{n\geq 2} A_nu_n^2, \quad \forall\; u\in C_c(\N), u_0 = 0,
\end{align}
where
\begin{align*}
A_n  = H_{-2}\left(\frac{1}{n}\right) = \frac{1}{n^2}\left[1 + \left(1 - \frac{1}{n}\right)^{-2} - \left(1 - \frac{1}{n}\right)^{-\frac{1}{2}} - \left(1 + \frac{1}{n}\right)^\frac{3}{2}\right] = \sum_{k \geq 2} \frac{a_k}{n^{k+2}},
\end{align*}
with
\begin{align*}
a_k & = (k+1) - \frac{\binom{2k}{k}}{4^k} - (-1)^k\frac{3\binom{2k-4}{k-2}}{4^{k-1}k(k-1)}> 0, \quad \forall \; k \geq 2.
\end{align*}
Denote $v = {\rm div}(u)$, so $v_0 = 0$. Applying \eqref{kpp1} to $v$, there holds
\begin{align*}
\sum_{n \geq 1} |\Delta u_n|^2  = \sum_{n \geq 1} |\nabla v_n|^2
& \geq \sum_{n \geq 1}\omega_n v_n^2\\
& = \sum_{n \geq 2}\omega_{n-1}|\nabla u_n|^2\\
& = \frac{1}{4}\sum_{n \geq 2} \frac{|\nabla u_n|^2}{(n-1)^2} + \sum_{k = 2}\sum_{n \geq 2}\frac{\binom{4k}{2k}}{(4k-1)2^{4k-1}(n-1)^{2k}}|\nabla u_n|^2\\
& \geq  \frac{1}{4}\sum_{n\geq 2} A_nu_n^2 + \sum_{k = 2}\sum_{n \geq 2}\frac{(2k+1)^2\binom{4k}{2k}}{(4k-1)2^{4k+1}}\frac{u_n^2}{n^{2k+2}}.
\end{align*}
We used \eqref{a=-2bis} and \eqref{a<0} for the second inequality. More precisely, we have
\begin{align}
\label{rellich-2}
\sum_{n \geq 1} |\Delta u_n|^2 \geq \sum_{n \geq 2} \left(\frac{9}{16n^4} + \frac{15}{16 n^5} +  \frac{213}{128n^6} + \ldots\right)u_n^2, \quad \forall\; u\in C_c(\N), u_0 = u_1 = 0.
\end{align}
The above estimate could be still improved, an interesting question is to know whether the constant $\frac{15}{16}$ is sharp.

\medskip
A much more challenging question is to search optimal weights for $\ell = 2$. Seeing \eqref{rellich-f}, we could attempt to find suitable positive function $f \in C(\N)$ and use the corresponding weights $w = \Delta^2 f/f$, but we have no idea about the {\sl right} choice of $f$.

\medskip
In \cite{GKS}, with the factorization method, Gerhat-Krej\v{c}i\v{r}\'\i k-\v{S}tampach obtained the weights $\widetilde w_n$ in \eqref{gks} by choosing $\widetilde w = (\Delta^2 g)/g$ with $g_n = n^\frac{3}{2}$. However, their result is not a direct consequence of \eqref{rellich-f} (or Theorem \ref{thm2.1} with $V \equiv 1$ and $f = g$), because $\Delta g_n < 0$ for $n \geq 1$.


\subsection{} It's worthy to mention that situations in $\Z^d$ ($d \geq 2$) have different nature comparing to the continuous situation $\R^d$ for $d \geq 2$, especially when $d$ goes to infinity, even for the first order case. More precisely, consider the Hardy-Rellich inequalities
\begin{align}
\label{rellichZd}
\sum_{n\in \Z^d} |\Delta^\frac{\ell}{2}u_n|^2  \geq C_{\ell, d}\sum_{n\in \Z^d} \frac{u_n^2}{n^{2\ell}}, \quad \forall\; u \in C_c(\Z^d), u(0) = 0.
\end{align}
In \cite{Gup2}, Gupta proved recently that for any $\ell \geq 1$, the sharp constant in \eqref{rellichZd}, denoted by $C_{\ell, d}$, satisfies $C_{\ell, d} \sim d^\ell$ as $d \to \infty$, that is
$$0 < \liminf_{d\to\infty} \frac{C_{\ell, d}}{d^\ell} \leq \limsup_{d\to\infty} \frac{C_{\ell, d}}{d^\ell} < \infty.$$
Recall that in $\R^d$, we have the following sharp Hardy-Rellich inequalities (see \cite{OW})
\begin{align}
\label{rellichRd}
\int_{\R^d} |\Delta^\frac{\ell}{2}\varphi(x)|^2dx \geq 4^{-\ell}\prod_{k=1}^\ell (d - 2k)^2\int_{\R^d} \frac{\varphi(x)^2}{|x|^{2\ell}}dx, \quad \forall\; \varphi \in C_c^\infty(\R^d\backslash\{0\}).
\end{align}
Therefore, the sharp constant has the asymptotic behavior as $4^{-\ell}d^{2\ell}$ as $d$ goes to infinity.

\bigskip
\noindent
{\bf Acknowledgements.} The authors are partially supported by NSFC (No.~12271164) and Science and Technology Commission of Shanghai Municipality (No.~22DZ2229014). The authors would like to thank professor Bobo Hua for interesting discussion. The authors thank also the anonymous referees for their careful reading and valuable comments.

\end{document}